%% file: teichtri.tex
\documentclass{amsart}
\usepackage[dvips]{graphicx,color}

\newcommand{\ext}{\mathop{\rm ext}\nolimits}
\newcommand{\injrad}{\mathop{\rm injrad}}
\newcommand{\CC}{\ensuremath{\mathcal{C}}}
\renewcommand{\SS}{\ensuremath{\mathcal{S}}}

\newcommand{\hol}{\operatorname{hol}}
\newcommand{\uhol}{|\! \hol \! |}
\newcommand{\Mod}{\mathop{\rm Mod}}

\newcommand{\R}{\ensuremath{\mathbb{R}}}
\renewcommand{\H}{\ensuremath{\mathbb{H}}}
\newcommand{\Z}{\ensuremath{\mathbb{Z}}}
\newcommand{\N}{\ensuremath{\mathbb{N}}}

\newcommand{\PP}{\mathbf{P}}
\newcommand{\limni}{\lim\limits_{n\to\infty}}

\theoremstyle{remark}
\newtheorem*{remark}{Remark}
\newtheorem*{remarks}{Remarks}

\newtheorem{fact}{Fact}

\theoremstyle{definition}
\newtheorem*{definition}{Definition}
\newtheorem{theorem}{Theorem}
\newtheorem*{thma}{Theorem A}
\newtheorem*{thmb}{Theorem B}
\newtheorem{corollary}[theorem]{Corollary}
\newtheorem{proposition}[theorem]{Proposition}
\newtheorem{lemma}[theorem]{Lemma}

\begin{document}

\title[Thin triangles and a M.E.T.\ for Teichm\"{u}ller
geometry]{Thin triangles and a multiplicative ergodic theorem for
Teichm\"{u}ller geometry}
\author{Moon Duchin}

\maketitle

\section{Introduction}
\subsection{Overview}
In this paper, we prove a curvature-type result about
Teichm\"{u}ller space, in the style of synthetic
geometry.\footnote{{\em Synthetic geometry} encompasses metric
geometry as well as a range of axiomatic approaches to the study of
geometric objects, as opposed to the coordinatized or tensorial
approaches of classical differential geometry. Active areas in the
synthetic tradition include CAT$(0)$ spaces, $\delta$-hyperbolicity,
and geometric group theory in general.  For a foundational text of
modern synthetic geometry, see Busemann \cite{busemann}.} We show
that, in the Teichm\"{u}ller metric, ``thin-framed triangles are
thin"---that is, under suitable hypotheses, the variation of
geodesics obeys a hyperbolic-like inequality.  This theorem has
applications to the study of random walks on Teichm\"{u}ller space.
In particular, an application is worked out for the action of the
mapping class group:  we show that geodesics track random walks
sublinearly.

Recall that the Teichm\"{u}ller
space $T_{g,n}$ is a parameter space for marked metrics
on oriented surfaces of a fixed topological type
(the marking is a choice of generators for $\pi_1$; the type $(g,n)$
is the genus and number of punctures or boundary components, chosen so
that $\Sigma_{g,n}$ is a hyperbolic surface).
The points of $T_{g,n}$ are conformal classes of marked metrics---equivalently,
since there is exactly one metric of constant curvature in each conformal
class, each point can be identified with a (marked)
Poincar\'{e} metric on the surface $\Sigma_{g,n}$.
The mapping class group ${\rm Mod}(g,n)$
is the collection of isotopy classes of orientation-preserving diffeomorphisms of $\Sigma_{g,n}$.
For a small value $\epsilon$, the {\em cusps} of $T_{g,n}$ are the
regions containing metrics on $\Sigma_{g,n}$
with some nontrivial curve shorter than $\epsilon$; the complement of the cusps is called
the {\em thick part}.

For any two points of Teichm\"{u}ller space, there are quasiconformal maps between them;
Teichm\"{u}ller showed that there is a unique quasiconformal map of minimal dilatation
(the eccentricity of its ellipse field) ~\cite{teich}.  He defined a distance function accordingly,
though it yields only a Finsler---not a Riemannian---metric.
This Teichm\"{u}ller metric is one of several natural metrics on $T_{g,n}$.
Here, we will restrict attention to this choice of metric, and to
Teichm\"{u}ller spaces $T_g$ of compact hyperbolic surfaces.

There is a long and involved history of studying the metric geometry of Teichm\"{u}ller space
since the introduction of the Teichm\"{u}ller metric in the late 1930s.
In 1959, Kravetz ~\cite{kravetz} had a much-cited result purporting to show that
the Teichm\"{u}ller metric was
Busemann non-positively curved (that is, that it had a convex distance function).
Linch showed that this argument was
incorrect in her thesis of 1971 ~\cite{linch}, and Masur proved that Teichm\"{u}ller space was in fact not
Busemann non-positively curved and not CAT$(0)$
(a slightly stronger condition)
in his own thesis of 1975 ~\cite{m-thesis}.

\newpage
Since then, Teichm\"{u}ller space has been shown to have quite a few properties which are
suggestive of negative/non-positive curvature:

\begin{enumerate}
\item every finite group of isometries has a fixed point, as is the case for CAT$(0)$ spaces.
For $T_g$, the fixed point property
is equivalent to the Nielsen realization
problem and was settled affirmatively by Kerckhoff in 1983 ~\cite{k-nielsen}

\item $T_g$ admits a boundary at infinity similar to hyperbolic space---namely, the
visual and Thurston compactifications each give a sphere at infinity.
Kerckhoff showed those two to be slightly different ~\cite{k-thesis},
and the boundary theory has since been extensively developed by Masur and Kaimanovich-Masur
~\cite{masur-twoboundaries},\cite{k-masur},\cite{k-masur2}

\item the geodesic flow on moduli space (the quotient of Teichm\"{u}ller space by its isometry group,
the mapping class group) is ergodic ~\cite{masur-1982}, as is the case with quotients of
hyperbolic space; and geodesics in moduli space obey a logarithm law ~\cite{masur-loglaw}
governing their rate of escape which
is a direct analog of Sullivan's logarithm law for geodesics on hyperbolic manifolds
~\cite{sullivan-loglaw}

\item Teichm\"{u}ller space is relatively hyperbolic with respect to the cusps
(that is, the electric
Teichm\"{u}ller space -- obtained from $T_g$ by coning off at the cusps --
is $\delta$-hyperbolic) ~\cite{farb-thesis},\cite{mm1}
\end{enumerate}

However, other features
undercut the usual ways of asserting negative/non-positive curvature:
\begin{enumerate}
\item there are families of geodesic rays from every point
of $T_g$ which stay a bounded distance apart, so $T_g$ is not CAT$(0)$~\cite{m-thesis}
\item there are arbitrarily large geodesic triangles
for which one edge stays far from the opposing vertex, so $T_g$ is not $\delta$-hyperbolic
~\cite{masur-wolf}
\item there is sup behavior in the cusps:  the Teichm\"{u}ller metric in the cusps differs only
by an additive constant from the sup metric on a product of lower-dimensional Teichm\"{u}ller
metrics and some hyperbolic metrics
~\cite{minsky-product}.
The sup metric ensures some positive-curvature characteristics, such as the existence of many
pairs of conjugate points.  (That is, geodesics are not unique.)
\end{enumerate}

This state of affairs---negative curvature in the thick part, positive curvature in the
cusps---calls for a result accounting for the geodesic interactions between those two
parts of the thick-thin decomposition, since ``most" geodesic rays
are missed by restricting attention to one or the other.\footnote{That is,
for a generic point and direction, the associated geodesic winds deeply in and out of the cusps.
This is a consequence of Masur's logarithm law, already mentioned.}
This paper provides a result in that direction (Theorem A), establishing
a metric comparison condition for Teichm\"{u}ller space, expressed in terms of
geodesic triangles, which is not restricted to
geodesics in the thick part.  (See \S 2.1.)

Using this result,
we deduce a multiplicative ergodic theorem for Teichm\"{u}ller space (Theorem B),
furthering the parallels between Teichm\"{u}ller space and the theory of symmetric spaces and
nonpositively-curved metric spaces.
This application can be interpreted as evidence that Theorem A
is the right kind of condition to explain some of
the non-positive curvature properties of Teichm\"{u}ller space.

\subsection{Ergodic theory and the theory of random walks}

The Oseledec Theorem (or Multiplicative Ergodic Theorem)
describes the asymptotic behavior of a broad class of ergodic random walks,
and in particular is often quoted for products of random matrices  ~\cite{oseledec}.
If $G=SL_n(\R)$, say, then we can sample matrices by a probability measure
$\mu$ on $G$ and consider their product.  A random walk is formed on
the associated symmetric space $G/K$ by applying the successive products to a basepoint.
The theorem provides, in part, that under mild conditions on $\mu$ (finite first logarithmic moment),
the rates of exponential growth of the eigenvalues of the product are almost everywhere constant.
In terms of these invariants, called Lyapunov exponents, the M.E.T.\ also contains
a convergence statement on eigenspaces for the random walk.
This powerful theorem, first proved in the 1960s, has been foundationally important in
modern dynamics.

The same result can be usefully restated in geometric language
in the case of $SL_n(\R)$, for instance: the theorem asserts
that almost every sample path (for a random walk satisfying the hypotheses)
leaves its basepoint with a speed which is an invariant of the dynamical system
$(G,\mu)$.
Then when that speed -- the first Lyapunov exponent -- is positive,
almost every sample path is sublinearly
approximated by constant-speed travel along a ray in a flat.

Karlsson and Margulis proved a version of the multiplicative ergodic theorem
for nonpositively curved spaces in ~\cite{k-margulis}.

\begin{theorem}[Karlsson-Margulis]
Let $(Y,d)$ be a uniformly convex, complete metric space that is Busemann non-positively curved.
Let $\Gamma$ be a semigroup of semicontractions $D\to D$, where $D$ is
a nonempty subset of $Y$, and fix a point $y\in D$.
Let $(\Omega,\PP)$ be a measure space with $\PP(\Omega)=1$ and
let $L:\Omega\to \Omega$ be an ergodic,
measure-preserving transformation.
For a measurable map $\pi: \Omega\to \Gamma$, form the cocycle
$$u(n,\omega):= \pi(\omega)\cdot \pi(L\omega) \cdots \pi(L^{n-1}\omega),$$
and let $y_n=u(n)y$.
Then there is a value $A$ such that
for $\PP$-a.e.\ $\omega$, the limit $\limni \frac 1n d(y, y_n(\omega))$ equals $A$.
If $\PP$ has finite first moment
~{\rm(}$\int_\Omega d(y,y_1) \ d\PP < \infty${\rm)}~
 and $A>0$,
then for $\PP$-a.e.\ $\omega$, there is a unique geodesic $\gamma$
starting at $y$ and such that $$\limni \frac 1n d(y_n,\gamma(An)) =0.$$
\end{theorem}

In other words, consider a random walk by semicontractions
(distance non-increasing maps) on a space satisfying
the geometric hypotheses.  If the average jump size is not too large (finite first moment),
then there is a well-defined dominant speed of deviation from the basepoint (namely $A$).
When that speed is positive, it follows that
for almost every sample path, there is a unique geodesic ray starting at the
basepoint and deviating sublinearly from the path.

Below, we fashion a multiplicative
ergodic theorem for Teichm\"{u}ller space.
A metric comparison result (Theorem A) is proved and used in the proof
in place of the Busemann nonpositive
curvature assumption (which fails for Teichm\"{u}ller space, as noted above).
The mapping class group acts by isometries (which are, in particular, semicontractions),
and we consider a measure $\mu$ on $\Mod(g)$.

\begin{thmb}
Let $(Y,d)=(T_g,d_T)$ be Teichm\"{u}ller space with the Teichm\"{u}ller metric.
Let $K$ be some thick part of $T_g$ and fix a point $y\in K$.
Suppose $\mu$ is a probability measure on the mapping class group $\Gamma=\Mod(g)$
such that the group generated by its support is non-elementary (has no fixed points
on the boundary).
Let $\PP$ be the product measure on $\Omega=\Gamma^\Z$,
let $L:\Omega\to \Omega$ be the left-shift, and let $\pi(\omega)=\omega_0$ read off the
first component.
Form the cocycle
$$u(n,\omega):= \pi(\omega)\cdot \pi(L\omega) \cdots \pi(L^{n-1}\omega),$$
and let $y_n=u(n)y$.
Then there is a value $A$ such that
for $\PP$-a.e.\ $\omega$, the limit $\limni \frac 1n d(y, y_n(\omega))$ equals $A$.
If $\PP$ has finite first moment and $A>0$,
then for $\PP$-a.e.\ $\omega$,
there is a unique geodesic $\gamma$ starting at $y$ and such that
$$\lim_{n\to\infty}\left[\frac{1}{n}d(y_n,\gamma(An)) \cdot \chi_K(p_n) \right] = 0$$
for $p_n=\gamma(d(y,y_n))$.
\end{thmb}

This says, in other words, that
almost every sample path of the random walk by mapping classes
has an associated Teichm\"{u}ller
geodesic that tracks it sublinearly while it travels in the thick part.

In fact, as long as the geodesic $\gamma$ associated to $\{y_n\}$
leaves the thick part $K$ sublinearly, the $\chi_K$ term can be
dropped.\footnote{Note that in terms of visual measure on geodesics,
a logarithmic rate of escape is generic, so this might be a
reasonable hope. At the moment, however, the hitting measure $\nu$
of sample paths on the boundary has few properties known; whether it
is absolutely continuous with respect to visual measure is open.}

\subsection*{Acknowledgments}
Thanks especially to Alex Eskin and Howard Masur.

\section{Triangle comparison in the Teichm\"{u}ller metric}

\subsection{The ``thin-framed triangles are thin" condition}

\begin{figure}[ht]
\begin{center}
\input{variation.pstex_t}
\caption[Varying a point from a geodesic segment]{Can the point
$x$ be far from the geodesic segment $\overline{yz}$
while the distance $d(y,x)+d(x,z)$ remains close to $d(y,z)$?   \label{variation}}
\end{center}
\end{figure}
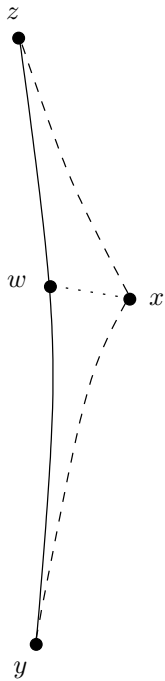

A natural way to think about the curvature of a metric space is to
ask about divergence of geodesics; flat spaces are characterized by
the linear spread of their geodesics, while positive and negative
curvature mean slower and faster divergence, respectively. A related
question one could pose of a fixed geodesic segment, say with
endpoints $y$ and $z$, would be to take a point $w$ on
$\overline{yz}$ and ask how much more distance is required to pass
from $y$ to $z$ through a perturbation $x$ of $w$ rather than
through $w$ itself.  (See Figure~\ref{variation}.) The property
proposed in this section will compare the extra distance needed
(namely $d(y,x)+d(x,z)-d(y,z)$) with the size of the perturbation
(namely $d(x,w)$).

Here, we will formulate this measurement in terms of triangles in
geodesic spaces.  Consider an arbitrary triangle $\triangle xyz$ in
a geodesic space $Y$. Relabeling the vertices if necessary, say that
the longest of the three pairwise distances between the points is
$d(y,z)$.  Let $w$ be the point along the geodesic from $y$ to $z$
such that $d(w,y)=d(x,y)$. For notational convenience, we will use
the letters $a,b,c,d$ for important distances: $c=d(y,z)$ is the
length of the longest side, $a=d(x,y)=d(w,y)$  and $b=d(x,z)$ are
the other two sidelengths, and $d=d(w,x)$ is the distance from $x$
to the specified point on the opposite side.

The following property was suggested by the
work of Karlsson and Margulis in \cite{k-margulis}.

\begin{definition} A collection $\mathcal{T}\subset Y^3$ of geodesic triangles in
the space $Y$ has the property that {\em thin-framed triangles are thin} if
\begin{center}
($\star$) there is some function $f(t)$ tending to zero as $t \to 0$
such that
$$\forall \rho>0, \quad a+b-c < a \cdot \rho \quad \implies \quad d < a \cdot f(\rho)$$
\end{center}
for every triangle in $\mathcal{T}$.
($\mathcal{T}$ will be taken to be all of $Y^3$ when not specified.)
\end{definition}

The idea behind using this definition is that the rate at which $f(t)$ approaches
zero (as $t$ goes to zero) detects the curvature of the space.

\newpage
\begin{proposition}
\begin{enumerate}
\item[\quad]
\item All triangles in trees (including both simplicial and $\R$-trees) satisfy ($\star$)
with bounding function $f(t)=t$, and no space has a bounding function going to zero
faster than linearly.
\item All triangles in $\R^2$ satisfy ($\star$) with bounding
function $f(t) = \sqrt{2t}$, and this function gives a sharp bound for the plane.
\item In $\delta$-hyperbolic spaces, if $\mathcal{T}$ is taken to be any collection
of triangles with sidelengths bounded away from zero, then ($\star$) is satisfied
with linear rate.
\item Generally, for a space $Y$,
$$Y ~\hbox{is}~ CAT(0)
\implies ~\hbox{thin-framed triangles are thin}~ \implies ~\hbox{geodesics are unique}.$$
In particular, spheres do not satisfy $(\star)$ for any bounding function.
\end{enumerate}
\end{proposition}

\begin{proof}
\begin{enumerate}
\item Graphs can be regarded as metric spaces with the path metric, assigning
distance one to each edge; non-degenerate triangles in trees are tripods (this holds
also for $\R$-trees).
In a tripod, $a+b-c=d$, so $f(t)=t$ is sharp.

In general, in any metric space, $c+2d$ is the length of one path from $y$ to $x$ to $z$
(going through $w$) and $a+b$ is the length of the most efficient path
from $y$ to $x$ to $z$, so
$d\ge \frac 12 (a+b-c)$ for any triangle in any metric space.
Thus a linear bound is best-possible (for any bounding function $f(t)$ for any metric space,
$\lim\limits_{t\to 0}\frac 1t f(t)$ is finite).

\item In a planar triangle, consider the angle $\theta$ at vertex $y$.  By the law
of cosines,
$$ \frac{2a^2-d^2}{2a^2} = \cos \theta = \frac{a^2+c^2-b^2}{2ac}$$

so that $$d^2 = \frac ac (b^2 - (c-a)^2) = a (a+b-c) \frac{c+b-a}c.$$
Now, using the hypothesis ($a+b-c<a\rho$)
and the fact that $b-a$ is less than $c$ but can get arbitrarily
close when $b\approx c$, we get that
$$d^2 < 2\rho a^2$$
and conclude that $d < \sqrt{2\rho} \cdot a$.

\item In any $\delta$-hyperbolic
space, such as the hyperbolic plane $\H^2$, the metric is only
boundedly far from that in a tree. For every nondegenerate triangle
in any metric space, there is an associated tree with three
endpoints (a {\em tripod}) whose edge lengths match the triangle.
One definition of $\delta$-hyperbolicity for a space $Y$ is that
$\delta$ gives a uniform bound for the {\em insize} of geodesic
triangles from $Y$ (the diameter of the set of three points in the
triangle which map to the central vertex in the tripod).

But this means that, if $r,s,t$ are chosen to be the lengths of the legs of the tripod
(that is, $a=r+s$, $b=s+t$, $c=r+t$),
then $d(w,x)\le 2s+\delta$.  Since $a+b-c=2s$, we have $d\le (a+b-c)+\delta$,
and under the hypothesis that $a+b-c<a\rho$, this yields $d\le a\rho+\delta$.

If we restrict attention to triangles with sidelengths bounded away from zero, say
$a\ge 1$, then $d\le a(\rho+\delta)$,
so the function $f(t)=t+\delta$ provides a bound as needed for ($\star$).

\item The definition of CAT$(0)$ provides that distances within a geodesic triangle are
less than or equal to those in a corresponding planar triangle.  So it is clear
that CAT$(0)$ spaces satisfy $(\star)$ with the same bound as we found for $\R^2$.

On the other hand, spaces with non-unique geodesics, and in particular spheres,
do not satisfy $(\star)$ for any bounding function:
construct a triangle by starting with a pair of
points that are joined by two distinct geodesics, and take the third
vertex to be any interior point on one of the geodesics (but not both).
For such a triangle, $c=b+a$, so
$\rho$ can be taken arbitrarily close to zero, but $d/a$ is non-zero.
\end{enumerate}
\end{proof}

\begin{remark} In the case of the hyperbolic plane $\H^2$,
we will not be able to improve on the $\R^2$ bound (where $f(t)\to 0$ at square-root rate)
overall, since
the metric on $\H^2$ degenerates to the Euclidean metric at small scale.
The proposition above, though, says that we can do better in the large:
large-scale triangles in negative curvature admit a linear bound for ($\star$).
This is different from the Euclidean case, where passing to the large scale does not
change the bounding function:  the square-root bound is still sharp.
\end{remark}

\medskip

Below (\S 2.3), a large-scale statement that thin-framed triangles
are thin will be proved for Teichm\"{u}ller space.

\subsection{Background}

General references for Teichm\"{u}ller theory, quadratic differentials, and polarized
flat structures can be found in \cite{abikoff}, \cite{gardiner},
\cite{strebel}, \cite{flp}/\cite{thurst-cpct}.

\medskip

One of many equivalent definitions of {\em Teichm\"{u}ller space} is
$$T_g = \left\{ (X,f) | \ f ~\hbox{quasiconformal}~ X_0\to X \right\}/\sim$$
where $X$ is a Riemann surface with genus $g$, $X_0$ is a fixed Riemann surface
serving as a basepoint,
and the equivalence relation is given by
$$(X_1,f_1)\sim(X_2,f_2) \iff \exists
~\hbox{conformal}~ h  ~\hbox{such that}~ h\simeq f_2\circ f_1^{-1}.$$
Accordingly, the Teichm\"{u}ller metric is given by
$$d_T((X_1,f_1),(X_2,f_2))= \frac 12 \ln \inf_{h} K(h),$$
infimized over quasiconformal $h \simeq f_2\circ f_1^{-1}$.
It is well-known that with this distance function, tangent spaces are equipped
with a notion of length (it is a Finsler metric),
but no notion of angle via an inner-product (it is not Riemannian).

The {\em mapping class group}, or modular group, is defined as
$${\rm Mod}(g) = \pi_0({\rm Diff}^+(\Sigma_g))={\rm Diff}^+(\Sigma_g)/{\rm Diff}_0(\Sigma_g).$$
The mapping class group acts on $T_g$ by changing markings and it is the full (orientation-preserving)
isometry group; the quotient ${\rm Mod}(g)\backslash T_g$  is the
{\em moduli space} $\mathcal{M}_g$ of (unmarked) conformal classes of metrics.  Because
the action has some fixed points, there is an orbifold structure on the quotient.

Among the metrics in each conformal class is a unique Poincar\'e metric (that is,
a metric of constant curvature $-1$) by uniformization.
Besides this metric, the conformal class contains many singular flat structures:  metrics
that are Euclidean away from a finite number of points, in which all the negative
curvature is concentrated.

A {\em quadratic differential} $\phi$ on a Riemann surface is a holomorphic 2-tensor
given in local coordinates by $\phi(z)dz^2$.
If the atlas for the Riemann surface is $(U_\alpha,z_\alpha)$, then the quadratic
differential is a system $(U_\alpha,z_\alpha,\phi_\alpha)$ of holomorphic functions
on each coordinate patch.  Given $\alpha$ and $\beta$, the functions transform by the rule
$\phi_\alpha(z_\alpha) \left(\frac{dz_\alpha}{dz_\beta}\right)^2 = \phi_\beta(z_\beta)$.
The collection of quadratic differentials
at a point $p\in T_g$ will be denoted $QD(p)$.  While the value of $\phi$ is not well-defined
at a point on the Riemann surface, its zeros (and their orders) are, because they are preserved by the
transformation rule.

It follows that the transition functions are semi-translations $(z\mapsto \pm z+c)$,
which means that the quadratic differential preserves a pair of line fields corresponding
to vertical and horizontal directions.  Note, however, that the directions are not oriented.
If instead we considered an abelian differential $f(z)dz$, we would obtain translation surfaces
with oriented vertical and horizontal directions.

Since a Euclidean coordinate can be obtained in each chart by $(\phi(z))^{1/2}$, a
quadratic differential gives rise to a {\em polarized flat structure} (the pair of line fields
is said to be a polarization)
with singularities at the zeros.  In fact, over an oriented topological surface of fixed
type,
there is a bijective correspondence between three sets of structures:
quadratic differentials, polarized flat structures, and Euclidean polygons with semi-translation
gluings and cone angles which are multiples of $\pi$.  Because of this correspondence,
the three notions will be used interchangeably.

Here, a flat structure will be written $X=(p,\phi)$ for some
quadratic differential $\phi\in QD(p)$, where $p\in T_g$.
The length measured in that flat metric may be denoted $\ell_X$,
$\ell_{(p,\phi)}$, or simply $\ell$, according to context.

A {\em saddle connection} is a curve on $X$ whose endpoints are at singularities and
which has no singularities in its interior.
A {\em cylinder} in a flat structure is an isometrically embedded Euclidean cylinder.
This occurs whenever there is a closed nonsingular curve because it must have
a family of homotopic parallel translates.  The boundary of a cylinder must consist of
a union of saddle connections.

For a closed curve $\gamma$ on a flat structure $X=(p,\phi)$, write
$\overline{\gamma}_X$ for its {\em straightening}---the geodesic in its homotopy
class---which
must take the form of a sequence of straight lines
$\gamma_1,\ldots,\gamma_k$.
Either $\gamma$ is part of a cylinder, in the case that
$\overline{\gamma}=\gamma_1$ is a single nonsingular straight-line
curve, or its geodesic representative is a sequence of saddle connections
$\gamma_1,\ldots,\gamma_k$ between not-necessarily-distinct singularities
(zeros of the quadratic differential $\phi$).

Where the flat structure $X=(p,\phi)$ is understood,
let $h_i$ and $v_i$ denote the components
of the {\em affine holonomy} of each segment:
$$\hol(\gamma_i) = (h_i,v_i) =({\rm Re} \int_{\gamma_i} \phi^{1/2} dz,
\ {\rm Im} \int_{\gamma_i} \phi^{1/2} dz).$$
Define the {\em unsigned holonomy} by
$|h|(\gamma)=\sum_i|h_i|$, \quad $|v|(\gamma)=\sum_i|v_i|$, and
$$\uhol(\gamma)=\left(|h|(\gamma),|v|(\gamma)\right).$$
A curve will be called {\em vertical} if $|h|=0$ and {\em horizontal}
if $|v|=0$.

\begin{remark} $|h|,|v|\le \left(|h|^2+|v|^2\right)^{1/2}\le \ell_X(\gamma)\le |h|+|v|$.
\end{remark}

\begin{fact} (Teichm\"{u}ller's Theorem)
Any two points of $T_g$ are connected by a unique geodesic.

The Teichm\"{u}ller geodesic flow is given by the multiplicative action of
$g_t=\left(\begin{smallmatrix}e^t&0\\ 0&e^{-t}\end{smallmatrix}\right)$ on quadratic differentials, so that
if a curve has piecewise holonomy coordinates $(h_i,v_i)$ in the flat structure given by
the pair $(x,\phi)$, then its piecewise coordinates in the push-forward by distance
$t$ are given by $(e^t h_i, \ e^{-t}v_i)$ on the flat structure
$g_t(X)=g_t(x,\phi)$ or $(x',g_t\phi)$ for the appropriate point $x'\in T_g$.
\end{fact}

That is, each direction in the visual sphere is associated to a quadratic differential,
and geodesic flow is given by the diagonal action of $SL_2(\R)$ on the corresponding flat
structure.  There is an appropriate normalization of quadratic differentials so that
the map of the open ball into Teichm\"{u}ller space, $B_1\hookrightarrow T_g$, is
a homeomorphism.  This is called the {\em visual embedding}.

\begin{remark} Unsigned holonomy behaves well under the Teichm\"{u}ller map:
$$\uhol_{g_t\bullet}=\left(e^t|h|,\ e^{-t}|v|\right).$$
\end{remark}

\medskip

We continue with key definitions.
For $x\in T_g$ and $\gamma\in\SS$ (the simple closed curves on $\Sigma_g$, up to homotopy),
the {\em extremal length} of $\gamma$ is
$$\ext_x(\gamma)=\sup_\ell \inf_{\gamma_0\sim\gamma}\ell(\gamma)^2,$$
where the sup is over metrics $\ell$ in the conformal class $x$.

\begin{fact} $\ext_x(\alpha)$ is achieved in a flat structure $X$ (called the Jenkins-Strebel differential)
which is composed
of a single cylinder, having $\alpha$ as the core curve.  This is discussed in \cite{strebel}.
\end{fact}

\begin{fact} Knowing how much curves are distorted is enough to compute Teichm\"{u}ller distance:
$$d_T(x,y) = \frac 12 \ln \sup_\alpha \left( \frac{\ext_x(\alpha)}{\ext_y(\alpha)} \right).$$
This formula is due to Kerckhoff \cite{k-thesis}.
\end{fact}

For a small value of $\epsilon>0$, the set $K=K(\epsilon)=\{x\in T_g : \injrad \ge \epsilon\}$
of metrics with no curve shorter than $\epsilon$
will be called a {\em thick part} of Teichm\"{u}ller space.
A thick part $K(\epsilon)$ projects to a
compact set $K_0=\pi(K)$ in moduli space but $K$ itself is not compact,
since it is ${\rm Mod}(g)$-invariant but
a sequence of mapping classes can send a point $x\in T_g$ to infinity by changing the marking.

\begin{fact}
For a thick part $K(\epsilon)$ of $T_g$, there is a constant
$F=F(\epsilon)$ (growing as $\epsilon \to 0$) such that $\forall \gamma \in \mathcal{S}$,
$$\frac 1F \ell_\psi(\gamma) \le \ell_\phi(\gamma) \le F \ell_\psi(\gamma)$$
for all $\phi, \psi \in QD(x)$.  In particular, this means that the length-squared of any
homotopy representative of $\gamma$ in any flat structure over the point $x$ is within
a factor $F^2$ of the extremal length $\ext_x(\gamma)$.
The proof is straightforward:  at each point in moduli space, there is some bound on the
greatest ratio of lengths of curves between metrics; this varies continuously over the
compact set $K_0$.
\label{F}
\end{fact}

\begin{lemma}\label{pantsplus}
There is some constant $R$ depending only on $\epsilon$ such that, for every point
$x$ in the thick part $K(\epsilon)$ of $T_g$, there is a collection of
closed curves of length $\le R$ whose complementary regions contain no
geodesic segment of length more than $3R$.
\end{lemma}

\begin{proof} Let
$$R_1=\sup_{x} \inf_{P\in\mathcal{P}} \max_{\gamma\in P} [\ext_x(\gamma)]^{1/2}~\quad
~\hbox{for}~x\in K=K(\epsilon)$$
$\gamma$ ranges over curves in the pants decomposition $P$, and $P$ ranges over
the pants complex $\mathcal{P}$.
The value $R_1$ exists because $\inf_{P\in\mathcal{P}} \max_{\gamma\in P} [\ext_x(\gamma)]^{1/2}$
varies continuously as $\pi(x)$ varies in the compact set $K_0$.
This means that there is some uniform constant ($F\cdot R_1$, for $F$ as in
Fact~\ref{F}) such that
every flat structure on every point over $K$ has some pants decomposition with
all pants curves shorter than that length.  Such a pants decomposition, $P$,
divides the surface into three-punctured spheres which,
in the Poincar\'{e} metric, are doubles of hyperbolic hexagons.  If the
$3g-3$ curves of $P$ are called {\em cuffs}, then these uniquely determine
$3g-3$ (not closed) curves called {\em inseams} which orthogonally connect
pairs of cuffs.  Call the cuffs $\gamma_i$ and the inseams $\delta_{ij}$
(for all $i,j$ such that $\gamma_i,\gamma_j$ bound the same pair of pants, including
$\delta_{ii}$ when $\gamma_i$ is separating).  In particular, the cuff lengths---here
no longer than $R_1$ and no shorter than $\epsilon$---determine the inseam lengths,
so let $R_2$ be largest possible inseam length.
Each $\delta_{ii}$ can be completed to a closed
loop of length at most $R_1+R_2$ by adding part of $\gamma_i$.
If there are indices $i_1,\ldots,i_k$ such that
$$\delta_{i_1i_2}, \ \delta_{i_2i_3}, \cdots, \ \delta_{i_ki_1}$$
is a sequence of inseams, then
there is a closed
loop of length at most $(3g-3)(R_1+R_2)$ containing those $k$ inseams as pieces.
Continuing until every inseam curve is included in some closed curve---a finite process---
let $R$ be the maximum of $R_1$ and any of the curves so obtained.
These closed curves, taken together with the pants curves, divide the
surface into simply connected pieces with six sides of length $\le R$, so no geodesic
segment on a complementary region can have length more than $3R$.
This implies that, if $\mathcal{M}$ is the multicurve given by the pants curves and
closed-up inseam curves, any geodesic $\gamma$ on $\Sigma$ satisfies the intersection
number inequality $i(\gamma,\mathcal{M})\ge \ell(\gamma)/3R$.
Finally, the geodesic representatives of these new curves are shorter or the same length,
and satisfy the same inequality.
\end{proof}

\begin{lemma}\label{worstcasebound}
For any two simple closed curves $\gamma_1,\gamma_2$ on $\Sigma_g$
and any point $x\in K(\epsilon)\subset T_g$,
the intersection number of the curves is bounded above with respect to the length in any metric
in the conformal class $x$:
$$i(\gamma_1,\gamma_2)\le \left(\frac{4}{\epsilon^2}\right)\ell(\gamma_1)\ell(\gamma_2).$$
\end{lemma}

\begin{proof}
Let $\gamma_{1}(j)$ be
a small piece of $\gamma_1$ having length less than $\epsilon/2$,
half the injectivity radius bound.  Between each two intersections
of $\gamma_2$ with $\gamma_{1}(j)$, the length along $\gamma_2$ is at least
$\epsilon/2$ (because otherwise a loop of length $<\epsilon$ would be formed).
Thus if $i(\gamma_{1}(j),\gamma_2)=n_j$, it follows that
$n_j\frac{\epsilon}2 \le \ell(\gamma_2)$.  Then cover $\gamma_1$ with
such pieces, of which at most $\frac{\ell(\gamma_1)}{\epsilon/2}$ are needed,
and let $n=\max\{n_j\}$.  This completes the proof:
$$i(\gamma_1,\gamma_2)\le n\cdot\frac{\ell(\gamma_1)}{\epsilon/2}\le
\left(\frac{4}{\epsilon^2}\right)\ell(\gamma_1)\ell(\gamma_2)$$
\end{proof}

\begin{fact}\label{wind} Generically, geodesics wind in and out of cusps.
This is discussed quantitatively in ~\cite{masur-loglaw}, where Masur shows that generically
(relative to the visual measure on directions)
a geodesic returns to the thick part between progressively deeper sojourns to various cusps.
He finds that $$\limsup_{t\to\infty}\frac{d_t}{\log t}=1/2,$$
where $d_t$ computed
by following the geodesic for time $t$, then measuring the distance in the moduli space
between the projections of the starting and ending points.
\end{fact}

\begin{fact} For a flat structure $X$,
the minimal components  under the flow $F_\theta$
(closed, flow-invariant sets on which the flow is minimal) are bounded by saddle connections
in the direction $\theta$.  Infinite leaves are dense in subsurfaces which are minimal
components for the flow, so it follows that infinite leaves get arbitrarily close
to saddle connections.  This is discussed in the Luminy lectures of Masur and Hubert-Schmidt
(\cite{hubert-schmidt},\cite{masur-luminy}).
\end{fact}

\begin{fact}For a flat structure $X$,
let $$V(X)=\{\mathbf{v}\in\R^2 ~\big|~ \exists ~\hbox{cylinder in}~ X
~\hbox{with holonomy}~ \mathbf{v} \}$$
and let $V_{\rm sc}(X)$ be defined similarly as the holonomies of saddle connections.
Then both are discrete as subsets of $\R^2$. \label{Vdiscrete}
The Luminy lectures (\cite{eskin},\cite{hubert-schmidt},\cite{masur-luminy}) contain a proof of this fact
and a discussion of related counting problems.
\end{fact}

\begin{fact}  Teichm\"{u}ller space is not $\delta$-hyperbolic. \label{masurwolf}
Masur and Wolf prove this in  \cite{masur-wolf}
by considering triangles given by a point $x\in T_g$ and its
Dehn twists about disjoint curves.  For large $n$, the triangle $\triangle x \tau_\alpha^n(x)
\tau_\beta^{-n}(x)$
does not have a uniform bound on the distance from the vertex $x$ to the opposite side.
(Ivanov later gave a simpler proof in \cite{ivanov}
using the observation that some geodesics spread slower than linearly,
which was the main finding of the much earlier paper \cite{m-thesis}.)
\end{fact}

\begin{fact} In a  flat structure $X$ with a minimal flow $F_\theta$, any
curve transverse to the flow (called a {\em section})
yields a {\em zippered rectangle} configuration by considering the first return of the flow
to the section (from either side):
this is a collection of rectangles with bases on the section
which recover the surface under gluings.
If the flat structure is a translation surface (the quadratic differential is the square
of some abelian differential; equivalently, the transition maps are translations) then
the flow is oriented and the rectangles are all distinct.  In general, the flow can
return to the same side of the section that it left from.
Zippered rectangles are a construction of Veech \cite{veech}.
\end{fact}

\subsection{A thin triangle theorem for Teichm\"{u}ller space}

The main theorem of this section states that, among
large-scale triangles in Teichm\"{u}ller space whose
four comparison points ($x,y,z,w$) are in some common thick part,
thin-framed triangles are thin.
In fact, the bounding function goes to zero not with the
$t^{1/2}$ rate found in the Euclidean case, but with the linear rate associated with negative curvature.

Consider a compact hyperbolic surface $\Sigma_g$
and points $x,y,z$ in $T_{g}$,
which all lie in some thick part $K(\epsilon)$, so that
there are no very short closed curves (shorter than $\epsilon$)
on the corresponding Riemann surfaces.
The Teichm\"{u}ller geodesics between such points, on the other hand, can wind deeply
in and out of cusps (Fact~\ref{wind})---that is, the shortest path between two
surfaces without short curves may certainly pass through surfaces with very short curves.

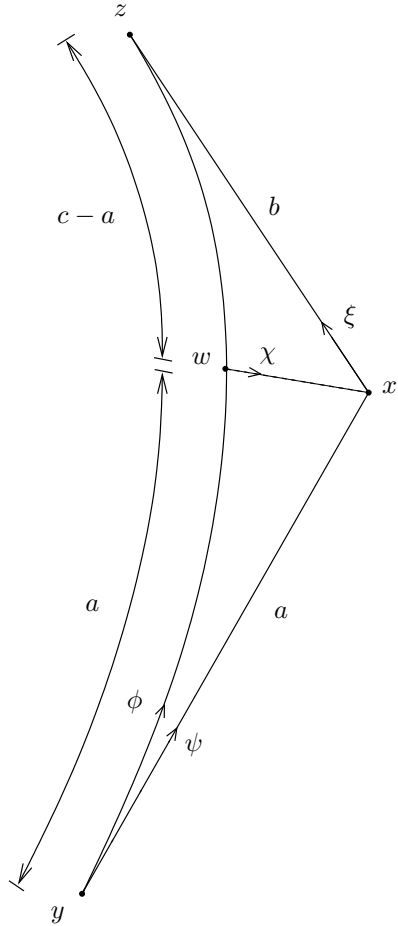
\begin{figure}[ht]
\begin{center}
\input{teichtri.pstex_t}
\caption{A triangle in Teichm\"{u}ller space. \label{maintri}}
\end{center}
\end{figure}

Choose quadratic
differentials $\phi,\psi$ (co-)tangent to the geodesics $\overline{yx}$ and $\overline{yz}$.
Then $(y,\phi)$ and $(y,\psi)$ are two polarized flat structures
with the same underlying Riemann surface.
See Figure~\ref{maintri}.

The goal for this section is

\begin{thma} Fix $K\subset T_g$, a thick part of Teichm\"{u}ller space.
Consider the collection of triangles
$\mathcal{T}_M=\{(x,y,z)\in K^3 | \ w\in K, \ a,b,c>M \}.$
For sufficiently large $M$,
$\mathcal{T}_M$ has the property that thin-framed triangles are thin,
satisfying ($\star$) with bounding function $f(t)=kt$, where $k$ depends on $K$.
\end{thma}

\begin{proof}[Proof of Theorem] We will consider a sequence of Teichm\"{u}ller geodesic
triangles $\{\triangle_i\}$ with sidelengths
$a=a(\triangle_i),b=b(\triangle_i),c=c(\triangle_i)$ as above.
We will assume that $a,b,c\to\infty$
and $\frac{a+b-c}a \to 0$
as $i\to\infty$.

To get an estimate of $d$, we define two families of curves for use in the analysis of
the flat structures of interest in the triangle.  We will consider
moderate-length curves at the vertices $y,z$.
An earlier lemma allows us to choose a length
cutoff which is large enough to assure that the families contain enough curves
to cut up the surface into simply connected pieces.

\begin{definition}
Where $\SS$ is the collection of simple closed curves on the surface,
up to homotopy, and $R$ is as in Lemma~\ref{pantsplus} above,
let

$$\CC_1 \ := \{\alpha\in\SS : \ell_\phi(\alpha) < R,  \ |v|_\phi/|h|_\phi <R \}$$

$$\CC_2 \ := \{\beta\in\SS : \ell_{g_c\phi}(\beta)<R, \ |h|_{g_c\phi}/|v|_{g_c\phi}<R \}$$
\end{definition}

\medskip

\begin{remarks}
\begin{itemize} \item[\ ]
\item These sets are nonempty; indeed, each contains a pants decomposition as
well as a collection of inseam-curves, as seen in the proof of Lemma~\ref{pantsplus}.
\item These sets are finite:  there are only finitely many homotopy classes of curves
such that any representative has length at most $R$ on a given flat structure $X$, by the
discreteness of the sets $V(X)$ and $V_{\rm sc}(X)$ (Fact~\ref{Vdiscrete}).
\item The slope bounds will be used to control the way the length changes under
Teichm\"{u}ller geodesic flow in the following lemma.
\end{itemize}
\end{remarks}

\medskip

We want to show $d/a\to 0$ as $\rho\to 0$.  Assume for contradiction that $d\gg a\rho$
for some $\rho$ such that $a+b-c<a\rho$.
In the following sequence of lemmas, we show that curves from
the families $\CC_1$ and $\CC_2$ are not changed too greatly in length between metrics based at
$w$ or at $x$ (Lemma~\ref{adeltabound});
that they are long and nearly perpendicular in the flat structure $(w,g_a\phi)$,
so have high intersection number (Lemma~\ref{lower});
that they would need to be nearly parallel
in $(w,\chi)$ in order to have $d\gg a\rho$ (Lemma~\ref{slopebounds});
and finally, that this would imply low
intersection number (Lemmas~\ref{tallrectangles} and ~\ref{upper}), producing a contradiction.

\begin{lemma}\label{adeltabound}
There is a constant $k_5=k_5(\epsilon)$
so that for all curves $\gamma$ in $\CC_1 \cup \CC_2$,
$$\frac 1{k_5}\cdot e^{-a\rho} < \frac{\ext_w(\gamma)}{\ext_x(\gamma)} < k_5
\cdot e^{a\rho}.$$
\end{lemma}

\begin{proof} Suppose $\alpha\in\CC_1$ has unsigned holonomy coordinates
$\uhol_{(y,\phi)}(\alpha)=(L,L')$,
where $L'/L<R$ by definition of $\CC_1$.

 Then, after applying geodesic flow for time $a$, we have
 $\uhol_{(w,g_a\phi)}(\alpha)=(e^aL,e^{-a}L')$.
The idea of the rest of this calculation is straightforward:
knowing that the length of $\alpha$
measured at $w$ is on the order of $e^aL$, we will
find the largest and smallest
possible length of $\alpha$ measured at $x$ by going around the triangle both ways from $y$.
On one hand, the distance from $y$ to $x$ is $a$, so the longest $\alpha$ could measure
at $x$ is $\sim e^a L$.  On the other hand, the length at $z$ is roughly $e^cL$ and the distance
from $z$ to $x$ is $b$, so the shortest possible length at $x$ is on the order of $e^{c-b}L$.
But then the biggest ratio of lengths is $e^{a+b-c}$; then the
``thin frame" hypothesis and Kerckhoff's formula will complete the proof.

 Now we will formalize this idea.  For large enough $a$,
 by repeated use of the
 observation that the length in any flat structure is bounded by $|h|\le \ell \le |h|+|v|$,
 we deduce the following collection of equalities and inequalities.

\begin{eqnarray*}
 e^aL \le \ell_{(w,g_a\phi)}(\alpha)\le 2e^aL,\\
 (e^aL)^2\le \ext_w(\alpha)\le (2Fe^aL)^2,\\
 \uhol_{(z,g_c\phi)}(\alpha)=(e^cL,e^{-c}L'),\\
 F^{-1} e^cL \le \ell_{(z,g_b\xi)}(\alpha) \le 2Fe^cL,\\
 F^{-1} e^{c-b}L \le \ell_{(x,\xi)}(\alpha),\\
 \left( F^{-1} e^{c-b}L \right)^2\le \ext_x(\alpha),\\
 \ell_{(y,\phi)}(\alpha)\le L+L',\\
 \ell_{(y,\psi)}(\alpha)\le F(L+L'),\\
 \ell_{(x,g_a\psi)}(\alpha)\le e^aF(L+L')\le e^aF(R+1)L,\\
 \ext_x(\alpha)\le \left( e^aF^2(R+1)L \right)^2.
\end{eqnarray*}

 Thus the length of $\alpha$ is not expanded by more than a constant factor
 when passing from $(w,\chi)$ to $(x,g_d\chi)$ and it may be contracted:
 $$ F^{-4}(R+1)^{-2} \le \frac{\ext_w(\alpha)}{\ext_x(\alpha)} \le
 \left( 2F^2 e^{a+b-c}\right)^2 < \left(2F^2 e^{a\rho}\right)^2.$$

 \medskip

In an exactly similar way, for $\beta \in \CC_2$,
$$\left(F^{-2}(R+1)^{-1}e^{-a\rho}\right)^2<  \left(F^{-2}(R+1)^{-1} e^{c-a-b}\right)^2
\le \frac{\ext_w(\beta)}{\ext_x(\beta)} \le 4F^{-2}$$
and both constants, $R$ and $F$, depend only on the choice of the compact part of
moduli space (that is, on $\epsilon$).
\end{proof}

\medskip

\begin{lemma}\label{lower}
There is a constant $k_6$ (depending on the genus $g$ and on $\epsilon$) so
that for any $\alpha$ in $\CC_1$, there is some $\beta$ in $\CC_2$
such that $$i(\alpha,\beta)\ge k_6 \ext_w(\alpha)\ext_w(\beta).$$
\end{lemma}

\begin{proof} We will consider the flat structure $(z,g_c\phi)$, where the $\CC_1$ curves
are long and near-horizontal and the
$\CC_2$ curves have moderate length.  By Lemma~\ref{pantsplus}, there is a collection of
curves from $\CC_2$ that decompose the surface into simply connected pieces, and the curves
are frequently intersected by long geodesics.
Call this collection of simple closed curves $\mathcal{B}$.

Fix any curve $\alpha\in\CC_1$.
For sufficiently large $c$, we have $\ell_{(z,g_c\phi)}(\alpha)\sim e^cA$.

It follows that $i(\alpha,\mathcal{B})\ge e^c \frac A{3R}$.
But there are no more than $6g-6$ curves in $\mathcal{B}$, so $\exists \beta\in \mathcal{B}$ such that
$$i(\alpha,\beta)\ge \frac{e^c A}{(18g-18)R}.$$
However, $\ext_w(\beta)\le e^{c-a}R$ and $\ext_w(\alpha)\le e^a A$, so
$$i(\alpha,\beta)\ge \frac{e^c AR}{(18g-18)R^2} \ge \frac{1}{(18g-18)R^2} \ext_w(\alpha)\ext_w(\beta).$$
\end{proof}

\begin{lemma}\label{slopebounds}
There are constants $k_{7},{k_{7}}'$ (depending on $\epsilon$)
so that
for any $\gamma\in\CC_1\cup\CC_2$, if we write $\uhol_{(w,\chi)}(\gamma)=(|h|,|v|)$, then
$$k_{7}
\cdot e^{d-a\rho} <  \frac{|v|}{|h|}  < {k_{7}}'\cdot e^{d+a\rho}.$$
\end{lemma}

\begin{proof} Let $X=(w,\chi)$.
First, note that $|v|$ is greater than $|h|$; otherwise, $\ell_X(\gamma)\sim |h|$ and
$\ell_{g_dX}(\gamma)\sim e^d |h|$,
contradicting Lemma~\ref{adeltabound} since $d$ is assumed much greater than $a\rho$.
Thus, we can write
$\uhol_{X}(\gamma)=(|h|,e^m |h|)$ where $m>0$, so that $\ell_X(\gamma)\sim e^m |h|$.
Observe that $m<2d$:  if not, then $\ell_{g_dX}(\gamma)\sim e^{m-d}|h|$, so the ratio of the
lengths again contradicts Lemma~\ref{adeltabound}.
Therefore, we have $\ell_{g_dX}(\gamma)\sim e^d |h|$, so
$\ext_w(\gamma)/\ext_x(\gamma) \sim e^{m-d}$.
But $m-d<a\rho \implies m<d+a\rho$ and $d-m<a\rho \implies m>d-a\rho$.
\end{proof}

The previous lemma says that the curves from the two principal families are near-vertical
considered in the flat structure $X=(w,\chi)$; they have slope close to $e^d$.
(Which isn't uniformly bounded because Teichm\"{u}ller space is not $\delta$-hyperbolic---Fact
~\ref{masurwolf}---and the triangles that witness the failure of $\delta$-hyperbolicity can be
arbitrarily large.)

\medskip

For any flat structure $X$, the surface decomposes in the vertical direction into a certain
number of Euclidean cylinders and a certain number of minimal components.  Fix a minimal
component $X_0$ and a section $\Gamma:[0,1]\to X_0$
transverse to the vertical foliation on $X$ and parameterized by arclength.
Let $\Gamma^{(l)}:[0,l]\to X_0$ be a subsection of length $l$.
Let $h(l)$ be the shortest return time to the subsection and
$w(l)$ be the minimum length of a subinterval of  $[0,l]$
on which the first return map is continuous.
(That is, $h(l)$ is the smallest height of a rectangle
and $w(l)$ is the smallest width of a rectangle
in the zippered rectangle configuration over $\Gamma^{(l)}$.)

\begin{lemma}\label{tallrectangles}
If $X_0$ is a minimal component of $X$ for the vertical flow, and $H>0$ is chosen arbitrarily,
then for sufficiently small $l=l(H)$, the vertical leaves over $\Gamma^{(l)}$ neither hit
singularities nor return to $\Gamma^{(l)}$ for at least time $H$.
\end{lemma}

\begin{proof}
Let $T$ be the interval exchange given by the first return to $\Gamma$ of the
vertical foliation.

Consider the points $T(0),T^2(0),\ldots,T^K(0)$ where $K=\lceil H/h(1) \rceil$.
Then for $l_0$ less than the minimum of these values, the first rectangle over $[0,l_0]$
has height at least $H$.  (It is possible to choose such an $l_0$
because none of the values $T^i(0)$ can be
equal to $0$ in a minimal component-- this would imply a closed orbit.)
Its width is at most $w(l_0)$, so by choosing $l_1<w(l_0)$ we can
be certain that all of the return times to the $l_1$-subsection are at least $H$.

Finally, by choosing $l_2<w(l_1)$, we ensure that the singularities are not encountered
until past height $H$ over the $l_2$-subsection.
\end{proof}

In particular, this lemma allows the construction of a
zippered rectangle configuration on $X_0$ such that all rectangles are
as tall as we like.

\begin{lemma}\label{upper}
Fix a constant $H>0$ and a flat structure $X$ in a thick part
$K(\epsilon)$ of Teichm\"{u}ller space.
Let $C_1,\ldots,C_r$ be the maximal
vertical cylinders in $X$ and let $X_1,\ldots,X_s$ be the minimal components
for the vertical flow.
Let $\Gamma_i^{(l_i)}:[0,l_i]\to X_i$ be sections in each minimal component
chosen relative to $H$ as in the previous lemma.
Then there is a constant $M=M(H,\{C_i\},\{\Gamma_i\},\epsilon)$
and a constant $k_{9}$
such that for any two curves $\alpha$ and $\beta$ with slope $|v|/|h|> M$,
$$i(\alpha,\beta) < k_{9}\frac{\ell_X(\alpha)\ell_X(\beta)}{H}$$
\end{lemma}

\begin{proof}
The intersection number is the minimum over homotopy representatives of the curves.
To compute it, it suffices to consider the geodesic straightening of
$\alpha$ and $\beta$ in the flat structure $X$.
We will denote by $i_C(\alpha,\beta)$ the number of intersections of the curves in the
interiors or on the boundary of any of the cylinders $C_i$ and by $i_Z(\alpha,\beta)$ the
number of intersections in or on the boundary of any of the minimal components, so that
$i(\alpha,\beta)\le i_C(\alpha,\beta)+i_Z(\alpha,\beta)$.

In the Euclidean picture $C_1 \sqcup C_r \sqcup X_1 \sqcup X_s$, each curve
decomposes into a certain number of maximal connected line segments, which we will call
simply {\em segments}.  See Figure~\ref{segments}.

\begin{figure}[ht]
\begin{center}
\includegraphics{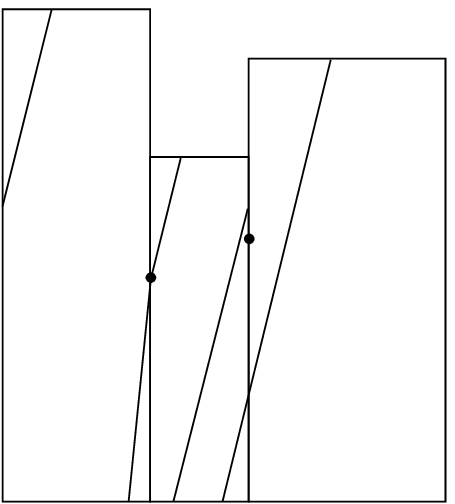} \quad \includegraphics{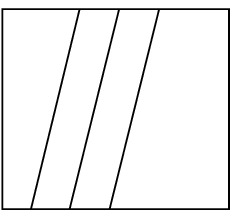}
\caption[A cylinder and a minimal component, with segments and sojourns of a curve indicated]{This figure
shows one cylinder, one minimal component made up of three rectangles,
and a total of eight segments of a single straightened curve.
Within the zippered rectangle figure,
sides of adjacent rectangles are glued together below the singularities.
The straightened curve can change slope only at singularities, as pictured.
The three segments that are shown in the cylinder are all part of a single sojourn
of the curve to the cylinder.\label{segments}}
\end{center}
\end{figure}

Let $w_0$ be the smallest of any of the widths of the rectangles or
cylinders; this depends only on the $C_i$ and $\Gamma_i$.  Let $M_1=H/w_0$.
We will show that if every segment of the curves $\alpha$ and $\beta$ has slope
$\ge M_1$, then $i(\alpha,\beta)< k_{\ref{upper}}\frac{\ell_X(\alpha)\ell_X(\beta)}{H}$.
In what follows, the length in the flat structure $X$ will be denoted simply $\ell$.

First note that each of the cylinders must be of (vertical) height at
least $\epsilon$, since the surface has no closed curves shorter than that.
Suppose $\alpha$ and $\beta$ each pass through cylinder $C_k$.
In fact, each can pass through the cylinder several times; let
$\alpha_1,\ldots,\alpha_p$ be the {\em sojourns} of $\alpha$ to $C_k$
(so that each is made up of many segments, one of which intersects the
left side of the cylinder and one of which intersects the right)
and likewise write
$\beta_1,\ldots,\beta_q$.  Each of these $pq$ curves has a fixed slope
for all of its segments
because there are no singularities in the cylinder, so let
$m_i$ be the slope of $\alpha_i$ and likewise $n_j$ of $\beta_j$
so that $m_i,n_j>M_1$ for all $i,j$.

Suppose $n_j>m_i$ for all $i,j$.
(If not, this estimate is similar but messier.)

Consider $i(\alpha_i,\beta_j)$.  Because $\beta_j$ has the greater
slope, at most two segments of $\alpha_i$ can intersect each of
its segments.  Thus $i(\alpha_i,\beta_j)\le 2 (w/h)n_j$, twice
the number of segments of $\beta_j$ in the cylinder with width $w$
and (vertical) height $h$.
But the larger the slope, the greater the length of $\beta_j$ must be:
$\ell(\beta_j)>n_jw$.  This means $i(\alpha_i,\beta_j)\le (2/h)\ell(\beta_j)$.
Repeating this calculation for
each of the pieces of $\alpha$, we find
$$i(\alpha,\beta_j)\le \frac{2p}{h} \ell(\beta_j)
\implies i(\alpha,\beta)\le \frac{2p}h \ell(\beta).$$
On the other hand, $\ell(\alpha_i)>|wm_i|>w_0M_1>H$, so $\ell(\alpha)>pH$,
and thus
$$i_{C}(\alpha,\beta)\le  \frac{2}{Hh} \ell(\alpha)\ell(\beta).$$

Let $A$ be the number of segments of $\alpha$ which
intersect the sections $\Gamma_i$
(that is, the bottom edges of rectangles in
zippered rectangle configurations).\footnote{Since the flow is not oriented, rectangles
may be repeated and the ``top" of one rectangle may be the ``bottom" of another.
But then treating them as distinct will if anything overcount the
number of intersections of the two curves.}
Since every such segment must have
length greater than $H$, it follows that $A<\ell(\alpha)/H$.

Every segment that intersects the top edge of a rectangle is naturally
paired (by the gluing pattern)
with one that intersects a bottom edge, so the number of such segments
is also $A$.

Finally, consider the segments which cross rectangles completely,
intersecting both sides.  The horizontal component is at least $w_0$,
so the vertical component---and therefore the length of the segment---is
at least $H$.  Thus there are no more than $A$ such segments.

Therefore, $\alpha$ has at most $3A$ Euclidean segments
in the zippered rectangles $X_1,\ldots,X_s$, and the same
bound holds for $\beta$.  Thus,
$$i_Z(\alpha,\beta)<  9A^2 < 9\frac{\ell(\alpha)\ell(\beta)}{H^2}.$$

To complete the proof, we need to account for segments whose slope may
be less than $M_1$.  We will choose a value of $M$ large enough that if
the curve $\alpha$ has overall
slope $M$, its segments of slope less than $M_1$ must be
very short.  Set $M>H(M_1+1)$.  Then segments of slope less than $M_1$ have
length $\ell_1\le |h|(M_1+1)$ (see Figure~\ref{bigslope}).

\begin{figure}[ht]
\begin{center}
\input{bigslope.pstex_t}
\caption[Very low-slope segments must be very short]{Having $|v|/|h|$ greater than $M$ means that
the pieces of slope less than $M_1$ are no more than $1/H$ of
the total length.
\label{bigslope}}
\end{center}
\end{figure}
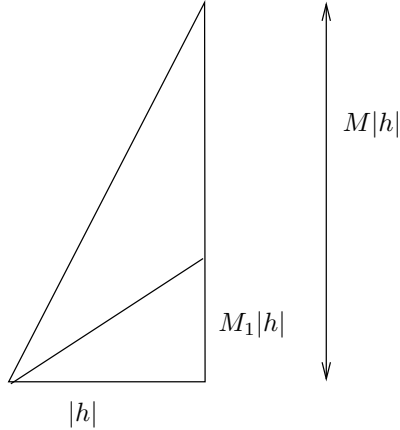

Since the entire curve has length
at least $M|h|$, this implies $\ell_1< \ell/H$.  But then, for the low-slope
pieces, Lemma~\ref{worstcasebound} implies that
$$i(\alpha,\beta)\le \left(\frac{4}{\epsilon^2}\right)\frac{\ell_X(\alpha)\ell_X(\beta)}{H}$$

This shows that the bound in the statement of the lemma holds for all cases,
completing the proof.
\end{proof}

Thus, by making $H$ (and thereby $M$) big enough, the intersection
number can be made less than an arbitrarily small proportion of
the product of the lengths of the curves.

\begin{remark}  This result is stable under perturbation; the same bound $M$ will
work for an open neighborhood of the flat structure $X$.
This is because a small perturbation (by say $\zeta$) will alter
the structures only by changing the widths of the rectangles and maybe
shifting the gluings (thus ruining the cylinders).  If these changes
are bounded by $\zeta$, which is taken to be very small relative to
the other parameters, the estimate is robust.
In particular, considering $i(\alpha_i,\beta_j)$ and supposing
$n_j>m_i$, the fact that a section of $\beta_j$ intersects at
most two sections of $\alpha_i$ is preserved as long as
$\zeta < \frac 12 \left( \frac{h}{n_j}-\frac{h}{m_i} \right)$.
But there are only finitely many values of the slopes, so $\zeta$ can
be chosen sufficiently small to preserve the estimates.
\end{remark}

\medskip

To finish the proof, we assemble the above information:
by Lemma~\ref{slopebounds}, curves
in $\CC_i$ would eventually have slope greater than any $M$.
But then Lemmas~\ref{upper} and ~\ref{lower} give contradictory inequalities
for the intersection numbers.
\end{proof}

Along the same lines, other triangle invariants can be bounded
relative to the framing defect $a+b-c$.  For instance, as a
``corollary of the proof," we can show that for large-scale
triangles, the insize is bounded above by $k(a+b-c)$ just as the
value $d$ is.  Such a bound (of the insize relative to the framing
defect) is obtained with a different formulation by Kent and
Leininger in a recent paper ~\cite{kent-lein}.

\section{Random Walk Application}

In this application, we consider random walks by mapping classes on Teichm\"{u}ller space.
The  triangle comparison property from the previous section is used to
argue that most sample paths are well-approximated by Teichm\"{u}ller geodesics, under
mild assumptions given below.

First, we use Teichm\"{u}ller-theoretic results of Masur and Kaimanovich-Masur to
associate to each (a.e.) sample path a geodesic that is a candidate to track the path.
Then, we invoke an ergodic argument from
Karlsson-Margulis to show that there are large, thin-framed triangles whose vertices
are points in the sample path.
Finally, we supply a geometric argument to complete the
estimate of the deviation of the sample path from its associated geodesic.

\subsection{Setup and ergodic ingredients}

The theorem of Karlsson and Margulis concerning ergodic transformations on nonpositively curved
spaces is stated in \cite{k-margulis}
in a great deal of generality.  We will present their results of interest
in full generality,
then indicate the chief and guiding example for this application, which will provide an opportunity
for concrete interpretation.

For their setup, let $(\Omega,\mu)$ be a measure space with $\mu(\Omega)=1$ and let $L:\Omega\to \Omega$
be a measure-preserving
transformation.  Let $a:\N\times \Omega\to \R$ be a subadditive (measurable) cocycle, that is
$$a(n+m,\omega)\le a(n,L^m \omega)+a(m,\omega)$$ for $n,m\in\N$,
$\omega\in \Omega$ (adopting the convention that $a(0,\omega)=0$).
We will assume that the following integrability condition is satisfied:
$$\int_\Omega a^+(1,\omega)\ d\mu(\omega)<\infty,$$
where $a^+(1,\omega)=\max\{0,a(1,\omega)\}$.  If $\mu$ satisfies this integrability
condition, then it is said to have {\em finite first moment}.
For each $n$, let $$a_n=\int_\Omega a(n,\omega)\ d\mu(\omega).$$
Then $A:= \lim_{n\to\infty} \frac 1n a_n$ exists and $A<\infty$.  (See ~\cite{k-margulis}.)

\begin{lemma}[Karlsson-Margulis ~\cite{k-margulis}, Prop 4.2]
\label{4.2}
Suppose that $L$ is ergodic and $A>-\infty$.
For any $\delta>0$, let $E_\delta$ be the set of $\omega\in \Omega$ such that
there exist an integer $N=N(\omega)$ and infinitely many $n$ such that
$$a(n,\omega)-a(n-k,L^k \omega) \ge (A-\delta)k$$
for all $k$, \ $N\le k\le n$.  If $E=\bigcup_{\delta>0}E_\delta$,
then $\mu(E)=1$.
\end{lemma}

\begin{corollary}[Kingmann Subadditive Ergodic Theorem, ~\cite{k-margulis} Cor 4.3]
\label{kingmann}
Under the same assumptions,
$$\lim_{n\to\infty} \frac 1n a(n,\omega)=A$$
for $\mu$-a.e. $\omega$.
\end{corollary}

A geometric interpretation is obtained by restricting attention to
the following case: Let $\Gamma={\rm Mod}(g)$ be the mapping class
group of a genus $g$ surface and let $\Omega$ be the state space
$\Omega=\Gamma^\Z$ so that an element $\omega\in\Omega$ is a
bi-infinite sequence of mapping classes $\omega=(\cdots, \omega_0,
\omega_1, \omega_2, \cdots)$. Let $\mathbf{P}$ be the product
measure on $\Omega$. Define the projection $\pi:\Omega\to \Gamma$ by
$\pi(\omega)=\omega_0$ and let the shift $L$ re-index by moving one
position to the left, i.e., $\pi(L^k\omega)=\omega_k$. (This is a
Bernoulli shift, so an ergodic transformation.) Let $Y=T_g$ be the
Teichm\"{u}ller space, fix a basepoint $y\in Y$, and let
$$y_n=y_n(\omega)=\pi(\omega)\pi(L\omega)\cdots\pi(L^{n-1}\omega)y
=\omega_0\omega_1\cdots \omega_{n-1}y.$$
This is the random walk being studied in this work:  starting at the basepoint $y$, the
subsequent points in the sample path are obtained by applying products of mapping classes to $y$.
Note the order of composition of the mapping classes is such that the first chosen is applied last.
As shorthand for this dynamical system, we will write the triple $(Y,\Gamma,\mu)$.

We will consider the cocycle $a(n,\omega)=d\bigl(y,y_n(\omega)\bigr)$, where $d$ is the
Teichm\"{u}ller distance.  Subadditivity is
an immediate consequence of the triangle inequality,
because $$a(n,L^m\omega)=d(y,\pi(L^m\omega)\cdots\pi(L^{m+n-1}\omega)y)=d(y_m,y_{m+n}).$$
The integrability assumption (finite first moment) for this distance cocycle is just the assumption
that the average jump size is finite.  Then $a_n$ can be interpreted as an average
distance from the basepoint after $n$ steps, so that $A$ should be thought of as a
dominant rate of escape for the random walk---an invariant of the dynamical system
$(Y,\Gamma,\mu)$.

\subsection{Teichm\"{u}ller ingredients}

Kaimanovich and Masur have showed that
under mild conditions on the measure $\mu$ (see below),
almost every sample path
$\{y_n\}$ of $(Y,\Gamma,\mu)$ converges to a uniquely ergodic foliation---call it $F$---on
the boundary of Teichm\"{u}ller space.
Thus, the unique geodesic from the basepoint $y$ to the limit foliation $F$ is a natural candidate to
approximate the sample path.

\begin{theorem}[Kaimanovich-Masur]
\label{candidateF}
If $\mu$ is a probability measure on the mapping class group $\Gamma=\Mod(g)$
such that the group generated by its support is non-elementary, then
for any
$y\in T_g$ and ${\mathbf P}$-a.e.\ sample path $\omega = \{\omega_n\}$ of the random walk $(\Gamma,\mu)$,
the sequence
$y_n$ converges in the Thurston compactification to a limit
$F = F(\omega) \in \mathcal{UE}\subset\mathcal{PMF}$.
\end{theorem}

This  originally appears as Kaimanovich-Masur~\cite{k-masur}, Thm 2.2.4, where more
is proved:  they show that this hitting measure $\nu$ on $\mathcal{PMF}$
(concentrated on $\mathcal{UE}\subset\mathcal{PMF}$, as above) is
the unique $\mu$-stationary probability measure $\nu$ on $\mathcal{PMF}$,
$\nu$ is purely non-atomic, and the measure space $(\mathcal{UE},\nu)$ is the Poisson
boundary of $(\Gamma,\mu)$.

\bigskip

We will also use a result of Masur.

\begin{theorem}\label{qdconv}
Suppose $F\in\mathcal{UE}\subset\mathcal{PMF}$ is a uniquely ergodic foliation.
Suppose further that there is a sequence of points $y_i\in T_g$ that converges to $F$
in the Thurston compactification.
Then, letting $y=y_0$ be a basepoint,
there is a sequence of quadratic differentials
$\varphi_i \in QD(y)$ so that for each $i$, $\exists t>0$ such that $y_i=g_t(y,\varphi_i)$,
and those converge in $QD(y)$ to a quadratic differential $\varphi$ such that
$F=\lim\limits_{t\to\infty}g_t(y,\varphi)$.
Therefore, for any fixed $m>0$, the points $g_m(y,\varphi_i)$ converge
to $g_m(y,\varphi)$ in $T_g$.
\end{theorem}

\bigskip

This follows from Masur ~\cite{masur-twoboundaries}, p.184,
a result comparing the visual and Thurston compactifications of $T_g$, which
were known by work of Kerckhoff not to be homeomorphic.  The boundary spheres
are the unit sphere of $QD$ and the sphere $\mathcal{PMF}$ of projective measured
foliations, respectively; Masur showed that $B_1 \cup Q_1^{ue} \cong T_g \cup \mathcal{UE}$,
where $Q_1^{ue}$ is the set of quadratic differentials in the unit sphere which have
uniquely ergodic vertical foliation.
In other words, the natural maps between $B_1$ and $T_g$ extend to the boundary spheres
at uniquely ergodic points.

\subsection{Assembling the ingredients}

Assume for this part that $A>0$ and fix an $\omega$ from the co-null subset of $\Omega$
such that
$\lim\limits_{n\to\infty} \frac 1n d(y,y_n(\omega))=A$.
Write $\gamma_k$ for the geodesic ray from $y$ through $y_k$ and $\gamma_\infty$ for the
geodesic ray from $y$ to $F=\lim y_n$.
For any $k$, let $r_k=d(y,y_k)$ and $p_k=\gamma_\infty(r_k)$.

\begin{lemma} \label{mainlemma}
For any $\delta>0$,  there exists $M_\delta$ with the property that
whenever $n$ is greater than $M_\delta$ there is some $m$
such that
\begin{eqnarray}
\label{thinframe}d(y,y_n)+d(y_n,y_m)-d(y,y_m)&\le&
\left( \frac{2\delta}{A-\delta} \right)d(y,y_n)\\
\label{bigenough}{\rm and} \ \ \ \ \ \ \ \ \ \ \ \ \ \  d(p_n,\gamma_m(r_n))&<& \delta.
\end{eqnarray}
\end{lemma}

\begin{proof}
We know that for any $\delta>0$,
there is an integer $N_\delta$ and infinitely many $m$ so that
$d(y,y_m)-d(y_n,y_m)\ge (A-\delta)n$ (Lemma~\ref{4.2})
and
$(A-\delta)n \le d(y,y_n) \le (A+\delta)n$ (Cor~\ref{kingmann})
for all $N_\delta\le n\le m$
(see Figure~\ref{ergodic}).

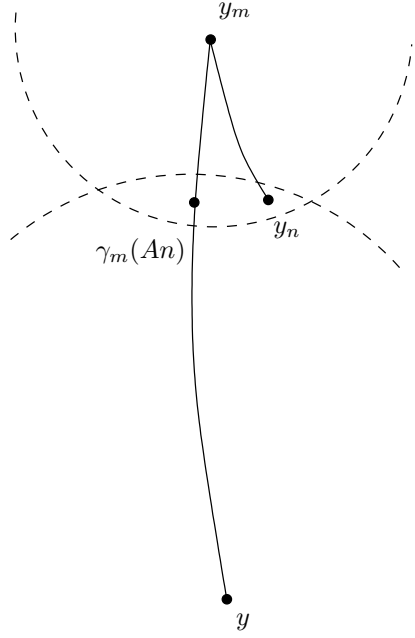
\begin{figure}[ht]
\begin{center}
\input{ergodic.pstex_t}
\caption[Geometric interpretation of inequalities used in the
ergodic argument]{The pair of inequalities above describes a region
around $\gamma_m(An)$ in which $y_n$ must fall.\label{ergodic}}
\end{center}
\end{figure}

Then $d(y,y_b)+d(y_n,y_m)-d(y,y_m)\le 2\delta n$.
Since $(A-\delta)n\le d(y,y_n)$, we finally get \eqref{thinframe},
which meets the ``thin frame"
hypothesis from the comparison triangle criterion discussed above
with respect to the triangle $\triangle yy_ny_m$.

For \eqref{bigenough}, recall that the points $\{y_i\}$
converge in the Thurston compactification to a foliation $F$ (Thm~\ref{candidateF}).
This means that the associated quadratic differentials converge in $QD$
(Thm~\ref{qdconv}). But by the visual embedding $B_1\to T_g$,
convergence of quadratic differentials
implies convergence in Teichm\"{u}ller space for points on the associated
geodesics at a fixed distance
from the basepoint.  This means that the points $\{\gamma_i(r_n)\}$ converge
to $\gamma_\infty(r_n)=p_n$;
choose sufficiently large $M_\delta\ge N_\delta$
so that $d(\gamma_m(r_n),p_n)<\delta$
for all $m>M_\delta$.

\end{proof}

\begin{lemma}\label{compactness} Suppose $\epsilon$ is fixed such that $K(\epsilon)$
contains $y$.
Then there is a constant
$k_{15}$ (depending only on $\epsilon$, $A$, and the genus $g$)
so that for any $\rho>0$ there exists $Q_\rho$ such that
$$n>Q_\rho, \quad p_n\in K \implies
\frac 1n d(y_n,\gamma_\infty(An))<k_{15}\cdot\rho.$$
\end{lemma}

\begin{proof}
Since $K=K(\epsilon)$ is the lift of a compact set in moduli space,
$y\in K$ implies that $y_i \in K$ for all $i$, since mapping classes
project to the identity on moduli space.

Choose $\delta$ small enough
that $\delta$, $\frac{2\delta}{A-\delta}$ are both less than $\rho$.

Choose an $\epsilon'$ small enough that $K(\epsilon')$ contains
a 1-neighborhood of $K(\epsilon)$.  (Thus $\epsilon'$ depends on $\epsilon$ and the genus.)
Then for $n$ bigger than $M_\delta$, there is an $m$ as in Lemma~\ref{mainlemma}
such that $d(p_n,\gamma_m(r_n))<\delta$, which is less than 1 without loss of generality.
For this value $m$,  it follows from $p_n\in K$ that $\gamma_m(r_n)\in K(\epsilon')$.

This means all four comparison points of $\triangle yy_ny_m$ are in some common thick part,
and $d(y,y_n)+d(y_n,y_m)-d(y,y_m)\le
\rho\cdot d(y,y_n)$ (by Lemma~\ref{mainlemma}\eqref{thinframe} and the choice of $\delta$).
We can then apply Theorem A, obtaining the estimate
$$d(y_n,\gamma_m(r_n))\le k_0\cdot \rho \cdot r_n,$$
where the constant $k_0$ depends on the value $\epsilon'$.

Putting this together, we find that
 $$d(y_n,\gamma_\infty(An)) \le d(y_n,\gamma_m(r_n)) + d(\gamma_m(r_n),p_n)
 + d(p_n,\gamma_\infty(An)) \le k_0\rho\ r_n + \rho + \rho n$$
$$\hbox{so}~ \frac{d(y_n,\gamma_\infty(An))}n \le \rho \left( \frac{k_0r_n}n + \frac 1n +1\right).$$

Set a new constant $k = 2(k_0A+1)$---this depends only on
the genus $g$ and the values $A$ and $\epsilon$.

Let $Q_\rho$ be the maximum of $M_\delta$ and a value large enough
that $\frac{k_0r_n}n + \frac 1n +1<k$ for all $n>Q_\rho$.  This is possible
because, once again applying Kingmann,
$\left(\frac{k_0r_n}n + \frac 1n +1\right) \longrightarrow (k_0A+1)$ as $n\to\infty$.
\end{proof}

\newpage
Together, these prove:

\begin{thmb}
Suppose $\mu$ is a probability measure on the mapping class group $\Gamma=\Mod(g)$
such that the group generated by its support is non-elementary.
Then if the measure has finite first moment and $A>0$,
it follows that for any
$y\in T_g$, any thick part $K$ of $T_g$ containing $y$,
and ${\mathbf P}$-a.e.\ $\omega$,

$$\lim_{n\to\infty}\left[\frac{1}{n}d(y_n,\gamma_\infty(An)) \cdot \chi_K(p_n) \right] = 0$$
\label{teich-met}
\end{thmb}

\end{document}

%% file: variation.pstex_t
\begin{picture}(0,0)%
\includegraphics{variation.pstex}%
\end{picture}%
\setlength{\unitlength}{4144sp}%
\begingroup\makeatletter\ifx\SetFigFont\undefined%
\gdef\SetFigFont#1#2#3#4#5{%
  \reset@font\fontsize{#1}{#2pt}%
  \fontfamily{#3}\fontseries{#4}\fontshape{#5}%
  \selectfont}%
\fi\endgroup%
\begin{picture}(855,4068)(1261,-4021)
\put(1306,-3976){\makebox(0,0)[lb]{\smash{\SetFigFont{10}{12.0}{\familydefault}{\mddefault}{\updefault}{\color[rgb]{0,0,0}$y$}%
}}}
\put(2116,-1771){\makebox(0,0)[lb]{\smash{\SetFigFont{10}{12.0}{\familydefault}{\mddefault}{\updefault}{\color[rgb]{0,0,0}$x$}%
}}}
\put(1268,-1666){\makebox(0,0)[lb]{\smash{\SetFigFont{10}{12.0}{\familydefault}{\mddefault}{\updefault}{\color[rgb]{0,0,0}$w$}%
}}}
\put(1261,-61){\makebox(0,0)[lb]{\smash{\SetFigFont{10}{12.0}{\familydefault}{\mddefault}{\updefault}{\color[rgb]{0,0,0}$z$}%
}}}
\end{picture}

%% file: teichtri.pstex_t
\begin{picture}(0,0)%
\includegraphics{teichtri.pstex}%
\end{picture}%
\setlength{\unitlength}{3947sp}%
\begingroup\makeatletter\ifx\SetFigFont\undefined%
\gdef\SetFigFont#1#2#3#4#5{%
  \reset@font\fontsize{#1}{#2pt}%
  \fontfamily{#3}\fontseries{#4}\fontshape{#5}%
  \selectfont}%
\fi\endgroup%
\begin{picture}(2356,5828)(430,-5961)
\put(2071,-1493){\makebox(0,0)[lb]{\smash{\SetFigFont{10}{12.0}{\familydefault}{\mddefault}{\updefault}{\color[rgb]{0,0,0}$b$}%
}}}
\put(923,-3983){\makebox(0,0)[lb]{\smash{\SetFigFont{10}{12.0}{\familydefault}{\mddefault}{\updefault}{\color[rgb]{0,0,0}$a$}%
}}}
\put(2108,-4028){\makebox(0,0)[lb]{\smash{\SetFigFont{10}{12.0}{\familydefault}{\mddefault}{\updefault}{\color[rgb]{0,0,0}$a$}%
}}}
\put(1181,-4591){\makebox(0,0)[lb]{\smash{\SetFigFont{10}{12.0}{\familydefault}{\mddefault}{\updefault}{\color[rgb]{0,0,0}$\phi$}%
}}}
\put(1546,-4861){\makebox(0,0)[lb]{\smash{\SetFigFont{10}{12.0}{\familydefault}{\mddefault}{\updefault}{\color[rgb]{0,0,0}$\psi$}%
}}}
\put(706,-5916){\makebox(0,0)[lb]{\smash{\SetFigFont{10}{12.0}{\familydefault}{\mddefault}{\updefault}{\color[rgb]{0,0,0}$y$}%
}}}
\put(1101,-241){\makebox(0,0)[lb]{\smash{\SetFigFont{10}{12.0}{\familydefault}{\mddefault}{\updefault}{\color[rgb]{0,0,0}$z$}%
}}}
\put(2546,-2176){\makebox(0,0)[lb]{\smash{\SetFigFont{10}{12.0}{\familydefault}{\mddefault}{\updefault}{\color[rgb]{0,0,0}$\xi$}%
}}}
\put(746,-1543){\makebox(0,0)[lb]{\smash{\SetFigFont{10}{12.0}{\familydefault}{\mddefault}{\updefault}{\color[rgb]{0,0,0}$c-a$}%
}}}
\put(1591,-2431){\makebox(0,0)[lb]{\smash{\SetFigFont{10}{12.0}{\familydefault}{\mddefault}{\updefault}{\color[rgb]{0,0,0}$w$}%
}}}
\put(2786,-2616){\makebox(0,0)[lb]{\smash{\SetFigFont{10}{12.0}{\familydefault}{\mddefault}{\updefault}{\color[rgb]{0,0,0}$x$}%
}}}
\put(2016,-2401){\makebox(0,0)[lb]{\smash{\SetFigFont{10}{12.0}{\familydefault}{\mddefault}{\updefault}{\color[rgb]{0,0,0}$\chi$}%
}}}
\end{picture}

%% file: bigslope.pstex_t
\begin{picture}(0,0)%
\includegraphics{bigslope.pstex}%
\end{picture}%
\setlength{\unitlength}{4144sp}%
\begingroup\makeatletter\ifx\SetFigFont\undefined%
\gdef\SetFigFont#1#2#3#4#5{%
  \reset@font\fontsize{#1}{#2pt}%
  \fontfamily{#3}\fontseries{#4}\fontshape{#5}%
  \selectfont}%
\fi\endgroup%
\begin{picture}(2010,2572)(939,-3706)
\put(1306,-3661){\makebox(0,0)[lb]{\smash{\SetFigFont{10}{12.0}{\familydefault}{\mddefault}{\updefault}{\color[rgb]{0,0,0}$|h|$}%
}}}
\put(2206,-3121){\makebox(0,0)[lb]{\smash{\SetFigFont{10}{12.0}{\familydefault}{\mddefault}{\updefault}{\color[rgb]{0,0,0}$M_1 |h|$}%
}}}
\put(2949,-1921){\makebox(0,0)[lb]{\smash{\SetFigFont{10}{12.0}{\familydefault}{\mddefault}{\updefault}{\color[rgb]{0,0,0}$M |h|$}%
}}}
\end{picture}

%% file: ergodic.pstex_t
\begin{picture}(0,0)%
\includegraphics{ergodic.pstex}%
\end{picture}%
\setlength{\unitlength}{4144sp}%
\begingroup\makeatletter\ifx\SetFigFont\undefined%
\gdef\SetFigFont#1#2#3#4#5{%
  \reset@font\fontsize{#1}{#2pt}%
  \fontfamily{#3}\fontseries{#4}\fontshape{#5}%
  \selectfont}%
\fi\endgroup%
\begin{picture}(2421,3837)(438,-3166)
\put(1801,-3121){\makebox(0,0)[lb]{\smash{\SetFigFont{10}{12.0}{\familydefault}{\mddefault}{\updefault}{\color[rgb]{0,0,0}$y$}%
}}}
\put(2026,-773){\makebox(0,0)[lb]{\smash{\SetFigFont{10}{12.0}{\familydefault}{\mddefault}{\updefault}{\color[rgb]{0,0,0}$y_n$}%
}}}
\put(1696,539){\makebox(0,0)[lb]{\smash{\SetFigFont{10}{12.0}{\familydefault}{\mddefault}{\updefault}{\color[rgb]{0,0,0}$y_m$}%
}}}
\put(962,-925){\makebox(0,0)[lb]{\smash{\SetFigFont{10}{12.0}{\familydefault}{\mddefault}{\updefault}{\color[rgb]{0,0,0}$\gamma_m(An)$}%
}}}
\end{picture}